\documentclass[12pt]{article}
\usepackage[a4paper, total={6in, 10in}]{geometry}
\usepackage{graphicx} 
\usepackage{soul}

\usepackage{amsmath, amsthm, amssymb, amsfonts, mathtools, mathrsfs}

\usepackage[utf8]{inputenc}
\usepackage{tikz-cd}
\usepackage[english]{babel}
\usepackage{etoolbox}
\usepackage{comment}
\usepackage{enumerate}
\usepackage{color}
\usepackage{bbm}
\usepackage{setspace}
\usepackage{bbm}
\usepackage[colorlinks]{hyperref}
\usepackage{float}
\usepackage{dsfont}
\usepackage{cancel}
\usepackage[nice]{nicefrac}
\usepackage{xfrac}
\usepackage[
backend=biber,
style=alphabetic,
sorting=ynt,
maxbibnames=99
]{biblatex}
\usepackage{xcolor}
\definecolor{Aranka}{RGB}{0, 135, 85}

\definecolor{Tomer}{RGB}{200, 20, 0}

\newtheorem{theorem}{Theorem}[section]
\newtheorem{lemma}[theorem]{Lemma}
\newtheorem{proposition}[theorem]{Proposition}
\newtheorem{corollary}[theorem]{Corollary}

\newtheorem{definition}[theorem]{Definition}
\theoremstyle{definition}

\newtheorem{remark}[theorem]{Remark}

\newcommand{\Acal}{\mathcal{A}}

\newcommand\NN{\mathbb N}

\newcommand\Hcal{\mathcal{H}}

\newcommand{\epsi}{\varepsilon}

\newcommand{\f}{\varphi}

\newcommand{\Prob}{\text{Prob}}

\newcommand{\ler}[1]{\left( #1 \right)} 
\newcommand{\lers}[1]{\left\{ #1 \right\}} 
\newcommand{\abs}[1]{\left| #1 \right|} 
\newcommand{\innprod}[1]{\left< #1 \right>} 

\newcommand{\acts}{\curvearrowright}

\newcommand{\dict}{\text{dict}}

\title{Minimal Models for Free Quasi-Invariant and Stationary Actions}
\author{Yair Hartman Tomer Zimhoni}
\date{November 2026}

\begin{document}
\maketitle

\begin{abstract}
    In this paper, we generalize the fact that essentially free p.m.p.\ actions of a countably infinite group admit a minimal model, proved by Weiss and Zucker. We prove that minimal models exist for every essentially free non-singular action. Furthermore, as in the p.m.p.\ case, we prove that all essentially free $\mu$-stationary actions (with respect to a fixed generating probability measure $\mu$) admit a common minimal model.
\end{abstract}

\section{Introduction}
    Let $G$ be a countably infinite group acting measurably on a standard Borel probability space $\ler{\Omega,\nu}$, and let $X$ be a compact metrizable $G$-space. We say that $X$ is a \textbf{model} of $\ler{\Omega,\nu}$ (or \textbf{models} $\ler{\Omega,\nu}$) if there exists a Borel probability measure $\eta \in \Prob\ler{X}$ such that $\ler{\Omega,\nu} \cong \ler{X,\eta}$ as measured $G$-systems.

    The discussion of ``universal'' minimal models for ergodic actions began with Weiss's result~\cite{weiss1989}. Weiss showed that every essentially free p.m.p.\ $\mathbb{Z}$-action admits a minimal model. Moreover, he constructed a single minimal, compact $\mathbb{Z}$-space, that is a model of all essentially free p.m.p.\ $\mathbb{Z}$-actions.

    Later on, in~\cite{weiss2012} Weiss extended this result much beyond $\mathbb{Z}$, and proved the same result for any countable infinite group.
    
    Extending Weiss's result Zucker~\cite{zucker2019note} presented a new construction in order to produce a single metrizable, minimal $G$-space that simultaneously models every essentially free p.m.p.\ $G$-system. He also showed that this construction is flexible enough to produce uncountable many different such models. This was generalized by Colin Jahel and Andy Zucker \cite{JahelZucker2021} to locally compact Polish groups.

    The first result of this paper extends~\cite{weiss2012} beyond the p.m.p.\ setting.
    
    \begin{theorem}\label{thm: any borel action has minimal model}
        Let $\ler{\Omega,\nu}$ be an essentially free, quasi-invariant measured $G$-system. Then there exists a metrizable minimal compact $G$-space $X$ which is a model of $\ler{\Omega,\nu}$.
    \end{theorem}

    Our second result generalizes~\cite{zucker2019note}, when we extend the settings of p.m.p.\ systems to that of $\mu$-stationary systems.

    \begin{theorem}\label{thm: any RW has universal minimal model}
        Let $\mu \in \Prob\ler{G}$ be a generating probability measure (as a semi-group). Then there exists a metrizable minimal compact $G$-space $X_\mu$ that is a model of every essentially free $\mu$-stationary $G$-system $\ler{\Omega,\nu}$.
    \end{theorem}

    \begin{remark}
        Since our subshift construction in Section~3 follows the mechanism of~\cite{zucker2019note}, varying the inductive parameters (specifically, the $\alpha$ sequences) yields $\aleph_1$-many mutually non-isomorphic minimal models. For clarity of exposition, we present the construction of a single minimal model in each case.
    \end{remark}

\subsection*{Acknowledgment}
This work was supported by the European Research Council (ERC)
under the European Union's Horizon Europe programme (BoundaryTheory, grant agreement No.~101078193).

\subsection*{AI declaration}
During the preparation of this work, the authors used ChatGPT (OpenAI) to assist with mathematical exploration, understanding Zucker's proof and verification of our results. The authors take full responsibility for all results, proofs, and references which appear in this paper.

\section{Preliminaries}

    Throughout this paper, let $G$ be a countable infinite discrete group with identity element $e$. 

    The main focus of this paper is the interplay between ergodic theory (measured $G$-systems) and topological dynamics (compact $G$-spaces). A measured $G$-system is a standard Borel probability space $\ler{\Omega,\nu}$ equipped with a Borel action $G \acts \Omega$. We say that $\nu$ (or the system $\ler{\Omega,\nu}$) is \textbf{quasi-invariant} if $g_*\nu \sim \nu$ for every $g \in G$, that is, $\nu(gA) = 0$ if and only if $\nu(A) = 0$ for every Borel set $A \subset \Omega$ and $g \in G$. A morphism of measured $G$-systems is a $G$-equivariant Borel map $f \colon \ler{\Omega_1,\nu_1} \to \ler{\Omega_2,\nu_2}$ with $f_*\nu_1 = \nu_2$. Such a morphism is an isomorphism if and only if it is essentially injective, namely, there exists a $G$-invariant $\nu_1$-conull Borel set $\Omega_0 \subset \Omega_1$ on which $f$ is strictly injective.

    A topological $G$-space is a compact metric space $X$ equipped with a continuous action $G \acts X$. Here morphisms are continuous $G$-equivariant maps.

    \begin{definition}
        Let $\ler{\Omega,\nu}$ be a $G$-system and let $X$ be a compact $G$-space. We say that $X$ is a \textbf{model} of $\ler{\Omega,\nu}$ if there exists a Borel probability measure $\eta \in \Prob\ler{X}$ such that
        \[
            \ler{\Omega,\nu} \cong \ler{X,\eta}
        \]
        as measured $G$-systems.
    \end{definition}

    \begin{remark}\label{remark: equivalent condition for a model}
        For a compact $G$-space $X$ to be a model of a $G$-system $\ler{\Omega,\nu}$, it suffices to find a $G$-equivariant Borel map $f \colon \Omega \to X$ that is essentially injective, in which case we set $\eta := f_*\nu$.
    \end{remark}

    Let $\mu \in \Prob(G)$ be a probability measure on $G$. We assume throughout that $\mu$ generates $G$ as semi-group, namely, every element of $G$ is a finite product of elements of $\operatorname{supp}\mu=\{g \mid \mu(g)>0\}$.  A Borel probability measure $\nu$ on a $G$-space $X$ is called \textbf{$\mu$-stationary} if
    \[
        \mu * \nu := \sum_{g \in G} \mu(g) g_*\nu = \nu.
    \]
    We denote the space of $\mu$-stationary probability measures on $X$ by $\Prob_\mu(X)$. Unlike $G$-invariant probability measures, which need not exist when $G$ is non-amenable, $\mu$-stationary measures always exist on any compact $G$-space $X$.

    In a similar way, a positive function $h:G \to \mathbb{R}_+$ is said to be \textbf{$\mu$-harmonic} if     
    \[
        h(\gamma) = \sum_{g \in G} \mu(g)h(\gamma g)
    \]
    We denote the collection of all non-negative $\mu$-harmonic functions by $\Hcal^{\ge 0}(G,\mu)$.

    Given a compact metric space $K$ (such as a finite alphabet $\Acal$ or the Cantor space $2^\NN$), the product space $K^G$ is compact under the product topology. We endow $K^G$ with the \textbf{right-shift action} of $G$, given for $g \in G$, $x \in K^G$, and $h \in G$ by
    \[
        (g \cdot x)(h) := x(hg).
    \]
    A closed, $G$-invariant subset $X \subset K^G$ is called a \textbf{subshift}. 

    When considering the alphabet $2^\NN$, a configuration $x \in (2^\NN)^G$ assigns to each group element $h \in G$ an infinite binary sequence $x(h) \in 2^\NN$. For each $n \in \NN$, we identify $2^n$ with the set of binary words of length $n$, and denote by $\pi_n \colon (2^\NN)^G \to (2^n)^G$ the coordinate projection that truncates each sequence $x(h)$ to its first $n$ entries:
    \[
        \pi_n(x)(h) := \ler{x(h)(1), \dots, x(h)(n)}.
    \]
    The map $\pi_n$ is continuous and $G$-equivariant, and commutes with restrictions to finite subsets $F \subset G$.

    Finally, we record two notions of separated sets that will be used to construct marker configurations. Let $\Omega$ be a Borel $G$-space, and let $L \subset G$ be a finite subset. A subset $B \subset \Omega$ is called \textbf{$L$-lacunary} if
    \[
        L b \cap L b^\prime = \emptyset \quad \text{for all distinct } b, b^\prime \in B.
    \]
    When $L$ is symmetric and contains $e$, this is equivalent to requiring that $b \notin L^2 b^\prime$ for all distinct $b, b^\prime \in B$. 

    Similarly, let $Q \subset G$ be a finite symmetric subset containing $e$. A subset $S \subset G$ is called \textbf{$Q$-spaced} if
    \[
        Q s \cap Q s^\prime = \emptyset \quad \text{for all distinct } s, s^\prime \in S,
    \]
    or equivalently, $s \notin Q^2 s^\prime$ for all distinct $s, s^\prime \in S$. A $Q$-spaced set $S$ is \textbf{maximal} (among $Q$-spaced sets with respect to inclusion) if and only if $Q^2 S = G$.

\section{Overview of Zucker's Proof}

    \subsection{Basic Models and Weak Models for Borel Actions}

    \subsubsection{The Space $\ler{2^\NN}^G$}
    We begin by describing the most elementary model of a $G$-system. Given a compact space $K$, the product space $K^G$ is a compact $G$-space endowed with the right shift action. Note that the space $K^G$ is never minimal, as it always contains constant functions (which are fixed points). 

    \begin{proposition}\label{prop: 2^N^G is always a model}
        The $G$-space $\ler{2^\NN}^G$ is a compact model of every $G$-system $\ler{\Omega,\nu}$.
    \end{proposition}

    \begin{proof}
        Without loss of generality, by passing to a Borel isomorphism, we may assume that $\Omega = 2^\NN$ (if in $\Omega$ is countable take an injection into $2^\NN$). Define
        \[
            U \colon \Omega \to \ler{2^\NN}^G, \quad U(\omega)(g) = g \cdot \omega.
        \]
        It is then straightforward to verify that $U$ is Borel, injective and $G$-equivariant, and therefore by Remark \ref{remark: equivalent condition for a model} $\ler{2^\NN}^G$ is indeed a model of $\ler{\Omega,\nu}$.
    \end{proof}

    Note that the fact that $\ler{2^\NN}^G$ is a model of a $G$-system $\ler{\Omega,\nu}$ depends neither on $\nu$ nor on the freeness of the action. 

    \subsubsection{Weak Models and the Space $Y_Q$}

    Although the following definition might appear ad hoc, it is an essential technical tool that we will use to construct minimal models.

    \begin{definition}
        Let $\ler{\Omega,\nu}$ be a $G$-system and let $X$ be a $G$-space. We say that $X$ is a \textbf{weak model} of $\ler{\Omega,\nu}$ if there exists a Borel $G$-equivariant map $\f \colon \Omega \to X$ (which is not required to be $\nu$-almost surely injective).
    \end{definition}

    Our next goal is to exhibit a family of weak models that will be useful in our context. Let $Q \subset G$ be a finite symmetric set. We say that $S \subset G$ is $Q$-spaced if for any distinct $s_1, s_2 \in S$, we have $Q s_1 \cap Q s_2 = \emptyset$. The proof of the following lemma is left to the reader. 

    \begin{lemma}
        A $Q$-spaced set $S$ is maximal (among $Q$-spaced sets) if and only if $Q^2 S = G$. Furthermore, for every finite symmetric set $Q\subset G$, maximal $Q$-spaced sets always exist.
    \end{lemma}
    
    We denote by $Y_Q$ the subshift of all maximal $Q$-spaced sets in $2^G$. Namely,
    \[
        Y_Q=\lers{s\in 2^G \mid s^{-1}\ler{1} \text{ is a maximal }Q\text{- spaced set}}.
    \]

    The following lemma will be used to prove the fact that $Y_Q$ is a weak model.

    \begin{lemma}\label{lemma: marker lemma}
        Let $\Omega$ be a free Borel $G$-space. Let $Q \subset G$ be a finite symmetric subset with $e \in Q$. Then there exists a Borel set $M \subset \Omega$ such that $Q m \cap Q m^\prime = \emptyset$ for all distinct $m, m^\prime \in M$, and $Q^2 M = \Omega$.
    \end{lemma}

    \begin{proof}
        Define a Borel graph on $\Omega$ by declaring $\omega, \omega^\prime$ to be adjacent if and only if $\omega = h \omega^\prime$ for some $h \in Q^2 \setminus \{e\}$. Since $Q^2 \setminus \{e\}$ is finite, this graph has bounded degree, and therefore admits a maximal independent Borel set $M$. The independence of $M$ implies that for any distinct $m, m^\prime \in M$, we have $m \notin Q^2 m^\prime$, which by the symmetry of $Q$ is equivalent to $Q m \cap Q m^\prime = \emptyset$. The maximality of $M$ ensures that $Q^2 M = \Omega$.
    \end{proof}

    \begin{corollary}\label{cor: Y_Q is weakly model}
        Let $\ler{\Omega,\nu}$ be an essentially free Borel $G$-system. Then for any finite symmetric subset $Q \subset G$ with $e \in Q$, $Y_Q$ is a weak model of $\ler{\Omega,\nu}$.
    \end{corollary}

    \begin{proof}
        By discarding a $G$-invariant null set if necessary, we may assume without loss of generality that the action of $G$ on $\Omega$ is strictly free. Let $M \subset \Omega$ be the Borel set provided by Lemma~\ref{lemma: marker lemma}. Define
        \[
            U_M \colon \Omega \to 2^G, \quad U_M(\omega)(g) = 1 \iff g\omega \in M.
        \]
        It is straightforward to verify that $U_M$ is Borel and $G$-equivariant. It remains to show that $U_M(\Omega) \subset Y_Q$. 
    
        Fix $\omega \in \Omega$, and consider the set $S = S_\omega := U_M(\omega)^{-1}(1)$. We show that $S$ is a maximal $Q$-spaced set. To see that $S$ is $Q$-spaced, let $s_1, s_2 \in S$ be distinct, and suppose for contradiction that $Q s_1 \cap Q s_2 \neq \emptyset$. Then there exist $q_1, q_2 \in Q$ such that $q_1 s_1 = q_2 s_2$. Since $s_1\omega, s_2\omega \in M$ and the action is free, we have $s_1\omega \neq s_2\omega$. However, $s_1\omega = q_1^{-1} q_2 s_2\omega \in Q^2 (s_2\omega)$, which contradicts the fact that $Q m \cap Q m^\prime = \emptyset$ for all distinct $m, m^\prime \in M$. Thus, $S$ is $Q$-spaced.

        Finally, we verify that $Q^2 S = G$. Given $g \in G$, the equality $Q^2 M = \Omega$ implies that $g\omega \in Q^2 M$. Hence, there exist $q_1, q_2 \in Q$ and $m \in M$ such that $g\omega = q_1 q_2 m$. Writing $m = (q_1 q_2)^{-1} g\omega \in M$, we see that $(q_1 q_2)^{-1} g \in S$. By the symmetry of $Q$, this yields $g \in Q^2 S$. 
    
        Therefore, $S_\omega$ is indeed a maximal $Q$-spaced set, establishing that $U_M(\Omega) \subset Y_Q$.
    \end{proof}

    Lastly, the following simple remark will be essential when constructing the space $Z_\bullet$ below.

    \begin{remark}\label{remark: product of models is model}
        If $X_1$ and $X_2$ are weak models of a $G$-system $\ler{\Omega,\nu}$, then the product $X_1 \times X_2$ is also a weak model of $\ler{\Omega,\nu}$. Moreover, if at least one of $X_1, X_2$ is a model, then $X_1 \times X_2$ is a model.
    \end{remark}

    \subsection{The Space $X_\bullet$} \label{subsection: the space X_alpha}

        \subsubsection{The Operator $\Phi$}
            We now focus our attention on subshifts of $\Acal^G$, where $\Acal$ is a finite alphabet. Let $X \subset \Acal^G$ be a subshift. For a finite set $F \subset G$, we denote by $\dict_X\ler{F} \subset \Acal^F$ the dictionary of $F$, which consists of all patterns supported on $F$. Formally,
            \[
                \dict_X\ler{F} = \lers{x|_F \mid x \in X}.
            \]
            
            For a finite set $D \subset G$, we say that $X$ is \textbf{$D$-irreducible} if whenever two finite sets $F_1, F_2 \subset G$ are $D$-apart (that is, $D F_1 \cap D F_2 = \emptyset$), any pair of patterns on $F_1$ and on $F_2$, appears simultaneously. Namely, for every $\alpha_1 \in \dict_X\ler{F_1}$ and $\alpha_2 \in \dict_X\ler{F_2}$, there exists some $x \in X$ such that $x|_{F_i} = \alpha_i$ for $i = 1, 2$.

            Let $F \subset G$ be a finite set. We say that $X$ is \textbf{$F$-minimal} if every possible pattern on $F$ appears in every element $x\in X$. Namely, for every $\alpha \in \dict_X\ler{F}$ and every $x \in X$, there exists some $g \in G$ such that $(g \cdot x)|_F = \alpha$. It is straightforward to verify that a subshift is minimal if and only if it is $F$-minimal for every finite $F \subset G$. 

            The next theorem provides a construction that, given a subshift, produces a new subshift that is ``more'' minimal while preserving the $F$-dictionary and maintaining irreducibility. The construction uses the subshift of maximal $C$-spaced configurations $Y_C$ described in Corollary~\ref{cor: Y_Q is weakly model}.

            \begin{theorem}[\cite{zucker2019note}, Section~3]\label{thm: existence of phi}
            Let $X \subset \Acal^G$ be a $D$-irreducible subshift for some finite symmetric subset $D \subset G$ with $e \in D$. Let $F \subset G$ be a finite symmetric set with $e \in F$, and let $E, C \subset G$ be sufficiently large finite symmetric sets with $e \in E, C$ such that $E^5 \subset C$ and $E$ contains at least $\abs{\dict_X\ler{F}}$ pairwise disjoint right translates of $D F$.
    
            Then there exists a continuous, $G$-equivariant map
            \[
                \f := \f\innprod{D,F,E,C} \colon X \times Y_C \to \Acal^G
            \]
            such that the subshift $\Phi = \Phi\innprod{X,D,F,E,C} := \f\ler{X \times Y_C} \subset \Acal^G$ satisfies:
            \begin{enumerate}
                \item $\Phi$ is $F$-minimal,
                \item $\dict_{\Phi}\ler{F} = \dict_X\ler{F}$, 
                \item every patern in $\dict_X\ler{F}$ occurs inside every pattern in $\dict_\Phi\ler{C^4}$,  
                and
                \item $\Phi$ is $C^5$-irreducible.
            \end{enumerate}
            Moreover, for every $s \in Y_C$, $x \in X$, and $g \in G \setminus E^3 s^{-1}\ler{1}$ (and in particular whenever $s^{-1}\ler{1} \cap E^4 g = \emptyset$), it holds that
            \begin{equation}\label{eq: phi preserving x out of s nbd}
                \f\ler{x,s}\ler{g} = x\ler{g}.
            \end{equation}
        \end{theorem}

        Intuitively, \eqref{eq: phi preserving x out of s nbd} asserts that each configuration $\f(x, s) \in \Phi$ agrees with $x \in X$ outside an $E^3$ right-neighborhood of the maximal $C$-spaced set $s^{-1}\ler{1}$.

    \subsubsection{The Construction of $X_\bullet$}

        We now construct a minimal subshift $X_\bullet \subset \ler{2^\NN}^G$, which serves as our candidate for a minimal model of a free Borel measured $G$-space.
        The data, denoted by $\bullet = (\textbf{D}, \textbf{E}, \textbf{F}, \alpha)$, is composed of three sequences of finite symmetric subsets of $G$, $\lers{D_n}, \lers{E_n}, \lers{F_n}$, and a sequence $\alpha \in \lers{0,1}^\NN$, satisfying the following conditions:
        \begin{itemize}
            \item $\abs{\alpha^{-1}\ler{1}} = \abs{\alpha^{-1}\ler{0}} = \infty$. We denote $T := \alpha^{-1}\ler{1}$ and $O := \alpha^{-1}\ler{0}$;
            \item $D_n$ contains the unit element and $D_n\subset E_n$
            \item $\bigcup_n D_n = \bigcup_n E_n = \bigcup_n F_n = G$;
            \item $E_n$ contains at least $2^{\abs{D_n}\ler{n+1}}$ pairwise disjoint right translates of $D_n^2$;
            \item $E_n^5 \subset F_n$; and
            \item $F_n^5 \subset D_{n+1}$.
        \end{itemize}

        We are now ready to define $X_\bullet$ inductively as a limit of subshifts $X_n \subset \ler{2^n}^G$. To do so, we embed $X_n$ into $(2^\NN)^G$ by padding configurations with zeros beyond the $n$-th coordinate.

        Let $X_0$ be the trivial subshift. For every $n \geq 1$, we define:
        \[
            X_n = 
        \begin{cases}
            X_{n-1} \times 2^G, & \text{if } \alpha\ler{n} = 0, \\
            \Phi\innprod{X_{n-1} \times 2^G, D_n, D_n, E_n, F_n}, & \text{if } \alpha\ler{n} = 1.
        \end{cases}
        \]
        We obtain that: 
    \begin{enumerate}
        \item For every $n$, $X_n$ is $D_{n+1}$-irreducible;
        \item Whenever $\alpha\ler{n} = 1$, $X_n$ is $D_n$-minimal; and
        \item Whenever $n < m$, $\dict_{\pi_n\ler{X_m}}\ler{D_{n+1}} = \dict_{X_n}\ler{D_{n+1}}$.
    \end{enumerate}

    Let $X_\bullet = \lim X_n$ be the resulting limit (In the Hausdorff topology). Item~1 ensures that the inductive construction of $X_n$ is well-defined, Item~3 implies that the sequence $(X_n)$ converges. Item~2, together with the infinitude of $T$ and the fact that every $D_n$ pattern occurs inside every $F_n^4$-pattern (which lie in $D_{n+1}$), ensures that $X_\bullet$ is minimal.

    \subsection{Zucker's Proof and Our Strategy}
        We now describe the essentials of Zucker's proof of the fact that for an appropriate choice of parameters, $X_\bullet$ is a global model for all essentially free p.m.p.\ actions of $G$. Throughout this subsection, $\bullet$ continues to denote the data $\bullet = \ler{\textbf{D},\textbf{E},\textbf{F},\alpha}$ fixed above. 

        We define the auxiliary $G$-space
        \[
            Z_\bullet = \prod_{n\in O} 2^G \times \prod_{n\in T} Y_{F_n} \cong \ler{2^\NN}^G \times \prod_{n\in T} Y_{F_n}.
        \]
        We denote a generic point in $Z_\bullet$ by $\ler{x,s}$, where $x = \ler{x_n}_{n\in O} \in \prod_{n\in O} 2^G \cong \ler{2^\NN}^G$ and $s = \ler{s_n}_{n\in T}$ with each $s_n \in Y_{F_n}$.

        The following lemma is a key technical tool of this paper. It appears in \cite{zucker2019note} as part of the proof of Proposition~8. For completeness and to adapt the statement to our setting, we provide a self-contained proof below.

        \begin{lemma}\label{lemma: key lemma}
    Let $\nu \in \Prob\ler{Z_\bullet}$. If for every $g \in G$ the following series converges:
    \[
        \sum_{n\in T} \nu\lers{\ler{x,s} \mid s_n^{-1}\ler{1} \cap E_n^4 g \ne \emptyset} < \infty,
    \]
    then $X_\bullet$ is a model of $\ler{Z_\bullet,\nu}$.
\end{lemma}

\begin{proof}
    Let
    \[
        Z_\bullet^\prime = \lers{\ler{x,s} \mid \forall g \in G, \; s_n^{-1}\ler{1} \cap E_n^4 g = \emptyset \text{ eventually}}.
    \]
    By the assumption in the lemma and the Borel--Cantelli lemma, $\nu\ler{Z_\bullet^\prime} = 1$. Therefore, it suffices to find an injective, $G$-equivariant Borel map $\psi \colon Z_\bullet^\prime \to X_\bullet$. Let $f_O \colon O \to O$ and $f_T \colon T \to T$ be two infinite-to-one surjections. Denote $\f_n = \f\innprod{D_n,E_n,F_n} \colon X_{n-1} \times 2^G \times Y_{F_n} \to X_n$. Consider the sequence of $G$-equivariant maps into $X_n$ 
    \[
        \begin{split}
            &\psi_n \colon Z_\bullet^\prime \to X_n, \\
            &\psi_n\ler{x,s} =
            \begin{cases}
                \ler{\psi_{n-1}\ler{x,s}, x_{f_O\ler{n}}}, & n \in O, \\
                \f_n\ler{\ler{\psi_{n-1}\ler{x,s}, s_{f_T\ler{n}}}, s_n}, & n \in T.
            \end{cases}
        \end{split}
    \]
    It is clear by definition that if $\psi_n$ converges, its image lies in $X_\bullet$. By the definition of $Z_\bullet^\prime$, for every $\ler{x,s} \in Z_\bullet^\prime$ and every $g \in G$ there exists some $n_0$ such that for all $n \geq n_0$, $s_n^{-1}\ler{1} \cap E_n^4 g = \emptyset$. Recall that by Theorem~\ref{thm: existence of phi}, whenever $s_n^{-1}\ler{1} \cap E_n^4 g = \emptyset$, we have $\f_n\ler{x,s_n}\ler{g} = x\ler{g}$. Hence, for every $n \geq n_0$,
    \begin{equation}\label{eq: psi_n is eventually insertion}
        \psi_n\ler{x,s}\ler{g} =
        \begin{cases}
            \ler{\psi_{n-1}\ler{x,s}, x_{f_O\ler{n}}}\ler{g}, & n \in O, \\
            \ler{\psi_{n-1}\ler{x,s}, s_{f_T\ler{n}}}\ler{g}, & n \in T.
        \end{cases}
    \end{equation}

    Thus, $\psi_n$ converges pointwise on $Z_\bullet^\prime$. We denote $\psi := \lim_{n\to\infty} \psi_n$.

    To show that $\psi$ is injective on $Z_\bullet^\prime$, consider a distinct pair $\ler{x,s}, \ler{x^\prime,s^\prime} \in Z_\bullet^\prime$. Suppose first that $s_k\ler{g} \ne s^\prime_k\ler{g}$ for some $k \in T$ and $g \in G$. Choose $n_0 \in \NN$ such that for all $n \geq n_0$, $s_n^{-1}\ler{1} \cap E_n^4 g = \emptyset$ and $\ler{s_n^\prime}^{-1}\ler{1} \cap E_n^4 g = \emptyset$. Since $f_T$ is an infinite-to-one surjection, there exists some $n_1 \geq n_0$ such that $f_T\ler{n_1} = k$. By \eqref{eq: psi_n is eventually insertion}, for any $n \geq n_1$, the configuration $\psi_n\ler{x,s}\ler{g} \in 2^n$ has its $n_1$-th entry equal to $s_k\ler{g}$, whereas the $n_1$-th entry of $\psi_n\ler{x^\prime,s^\prime}\ler{g}$ is $s^\prime_k\ler{g}$. To show that $\psi\ler{x,s} \neq \psi\ler{x^\prime,s^\prime}$ in the case where $x \ne x'$, one applies the exact same argument.
\end{proof}
        
    \subsubsection*{The Main Idea}
        Let $\ler{\Omega,\nu}$ be an essentially free $G$-system. By Proposition~\ref{prop: 2^N^G is always a model}, Corollary~\ref{cor: Y_Q is weakly model}, and Remark~\ref{remark: product of models is model}, we obtain that $Z_\bullet$ is a model of $\ler{\Omega,\nu}$. Write $f \colon \Omega \to Z_\bullet$ for the corresponding model map. Our strategy to prove Theorems~\ref{thm: any borel action has minimal model} and~\ref{thm: any RW has universal minimal model} is to choose parameters $\bullet = \ler{\textbf{D},\textbf{E},\textbf{F},\alpha}$ such that $\ler{Z_\bullet, f_*\nu}$ satisfies the conditions of Lemma~\ref{lemma: key lemma}. Consequently, considering the map $\psi \colon Z_\bullet^\prime \to X_\bullet$ from the proof above, the composition $\psi \circ f \colon \Omega \to X_\bullet$ is a $\nu$-essentially injective map; hence, $X_\bullet$ is indeed a model of $\ler{\Omega,\nu}$.

    \subsubsection*{Zucker's Proof Idea}

        In \cite{zucker2019note}, it is shown that if $\nu \in \Prob\ler{Z_\bullet}$ is a $G$-invariant measure, then for every $g \in G$ and $n \in \NN$,
        \[
            \nu\lers{\ler{x,s} \mid s_n^{-1}\ler{1} \cap E_n^4 g \ne \emptyset} \leq \frac{\abs{E_n^4}}{\abs{F_n}}.
        \]

        Now, we can add the condition that $\bullet$ satisfies
        \begin{equation}\label{eq: suff cond for pmp}
            \frac{\abs{E_n^4}}{\abs{F_n}} < \frac{1}{2^n}.
        \end{equation}
        Then the conditions of Lemma~\ref{lemma: key lemma} hold for every invariant measure $\nu \in \Prob\ler{Z_\bullet}$.

        In particular, let $\ler{\Omega,\nu}$ be a p.m.p.\ $G$-system and let $f \colon \Omega \to Z_\bullet$ be a model map. Since $f_*\nu$ is an invariant measure on $Z_\bullet$, we obtain that $X_\bullet$ is a model of $\ler{Z_\bullet, f_*\nu}$; composing $\psi$ with $f$ then shows that $X_\bullet$ is indeed a model of $\ler{\Omega,\nu}$.

\section{Minimal Model for any Free Borel Action}

    We are now ready to prove Theorem~\ref{thm: any borel action has minimal model}, asserting that any essentially free $G$-system admits a minimal model.

    \begin{lemma}\label{lemma: small shifted marker lemma}
        Let $G \acts \ler{\Omega,\nu}$ be a free Borel action, let $A, L \subset G$ be two finite subsets, and let $\epsi > 0$. Then there exist a finite, containing the unit element, symmetric subset $\Tilde{L} \subset G$ containing $L$ and a Borel subset $B \subset \Omega$ such that:
        \begin{enumerate}
            \item $B$ is $\Tilde{L}$-lacunary;
            \item $\Tilde{L}^2 B = \Omega$; and
            \item $\nu\ler{AB} < \epsi$.
        \end{enumerate}
    \end{lemma}

    \begin{proof}
        Take $N \in \NN$ such that $1/N < \epsi$, and choose $t_1, \dots, t_N \in G$ such that the sets $At_i$ are pairwise disjoint. Choose a finite, containing the unit element, symmetric set $L^\prime \subset G$ such that for every $i$,
        \[
            At_i, \, t_i^{-1}Lt_i \subset L^\prime.
        \]
        By Lemma~\ref{lemma: marker lemma}, we can choose a Borel set $M \subset \Omega$ that is $L^\prime$-lacunary and satisfies $\ler{L^\prime}^2 M = \Omega$.

        We now claim that the sets $At_i M$ are pairwise disjoint. Suppose that
        \[
            at_i m = a^\prime t_j m^\prime
        \]
        for some $a, a^\prime \in A$, $m, m^\prime \in M$, and $1 \leq i, j \leq N$. Since $at_i, a^\prime t_j \in L^\prime$, we have $L^\prime m \cap L^\prime m^\prime \ne \emptyset$, and thus by lacunarity, $m = m^\prime$. By the freeness of the action, we obtain $at_i = a^\prime t_j$, and since the sets $At_i$ are pairwise disjoint, it follows that $i = j$.

        Therefore,
        \[
            \sum_{i=1}^N \nu\ler{At_i M} \leq 1,
        \]
        and so for some $i$, $\nu\ler{At_i M} \leq 1/N < \epsi$. Setting
        \[
            B = t_i M \quad \text{and} \quad \Tilde{L} = t_i L^\prime t_i^{-1},
        \]
        we verify that all required properties hold:
        \begin{enumerate}
            \item $B$ is $\Tilde{L}$-lacunary: for every $b = t_i m \in B$ ($m \in M$), we have
            \[
                \Tilde{L}b = t_i L^\prime m,
            \]
            and the sets $t_i L^\prime m$ ($m \in M$) are pairwise disjoint;
            \item $\Tilde{L}^2 B = t_i \ler{L^\prime}^2 M = t_i \Omega = \Omega$; and
            \item $\nu\ler{AB} = \nu\ler{At_i M} < \epsi$.
        \end{enumerate}
    \end{proof}

    \begin{corollary}
        Theorem~\ref{thm: any borel action has minimal model} holds.
    \end{corollary}

    \begin{proof}
        Let $\ler{\Omega,\nu}$ be an essentially free Borel $G$-space, and let $K_n \nearrow G$ be an exhaustion of $G$ by finite sets. We inductively construct three sequences of finite symmetric subsets $D_n, E_n, F_n$ exhausting $G$ as described in Subsection~\ref{subsection: the space X_alpha}, with one additional requirement: at step $n$, we apply Lemma~\ref{lemma: small shifted marker lemma} with $L = E_n^5$, $A = \ler{E_n^4 K_n}^{-1}$, and $\epsi = \frac{1}{2^n}$ to obtain $F_n := \Tilde{L}$ and a corresponding Borel marker set $B_n \subset \Omega$ satisfying the assertions of the lemma.

        In this way, for every $g \in K_n$, we obtain
        \[
            \begin{split}
                {U_{B_n}}_{*}\nu
                \lers{s \in Y_{F_n} \mid 
                s^{-1}\ler{1} \cap E_n^4 g \ne \emptyset
                }
                &\leq 
                {U_{B_n}}_{*}\nu
                \lers{s \in Y_{F_n} \mid 
                s^{-1}\ler{1} \cap E_n^4 K_n \ne \emptyset
                }
                \\
                &= \nu \lers{\omega \in \Omega \mid \exists h \in E_n^4 K_n : h\omega \in B_n}\\
                &= \nu\ler{\ler{E_n^4 K_n}^{-1}B_n} < \frac{1}{2^n}.
            \end{split}
        \]
        Since every $g \in G$ eventually lies in $K_n$, it follows that
        \[
            \sum_{n \in T} {U_{B_n}}_{*}\nu
                \lers{s \in Y_{F_n} \mid 
                s^{-1}\ler{1} \cap E_n^4 g \ne \emptyset
                } < \infty.
        \]
        Thus, by Lemma~\ref{lemma: key lemma}, for this choice of parameters $\bullet$, $X_\bullet$ is a minimal model of $\ler{\Omega,\nu}$.
    \end{proof}

\section{Universal Minimal Model for Stationary Free Borel Actions}

    We now prove Theorem~\ref{thm: any RW has universal minimal model}, asserting that, given a generating random walk $\mu \in \Prob\ler{G}$, all essentially free $\mu$-stationary $G$-systems admit a common minimal model.

    \begin{lemma}\label{lemma: entries time is harmonic}
        Let $Q \subset G$ be a finite, symmetric subset with $e\in Q$, and let $\nu \in \Prob_\mu\ler{Y_Q}$. Set $p_\nu\ler{g} := \nu\lers{s \in Y_Q \mid s\ler{g} = 1}$. Then $p_\nu$ is $\mu$-harmonic and strictly positive.
    \end{lemma}

    \begin{proof}
        First, by the $\mu$-stationarity of $\nu$, we obtain harmonicity:
        \[
        \begin{split}
            p_\nu\ler{g} &= \nu\lers{s \in Y_Q \mid s\ler{g} = 1} \\
            &= \sum_{h \in G} \mu\ler{h} \nu\lers{s \in Y_Q \mid \ler{h \cdot s}\ler{g} = 1} \\
            &= \sum_{h \in G} \mu\ler{h} \nu\lers{s \in Y_Q \mid s\ler{gh} = 1} = 
            \sum_{h \in G} \mu\ler{h} p_\nu\ler{gh}.
        \end{split}
        \]
        Moreover, if $p_\nu\ler{e} = 0$, then since $\mu$ is generating, harmonicity and the non-negativity of $p_\nu$ imply that $p_\nu \equiv 0$. This means that $s^{-1}\ler{1} = \emptyset$ for $\nu$-almost every $s \in Y_Q$, which is impossible since the empty set is not a maximal $Q$-spaced set.
    \end{proof}

    We define $\Hcal^1(G,\mu) := \lers{h \in \Hcal^{\ge 0}(G,\mu) \mid h \geq 0, \, h\ler{e} = 1}$. In particular, by Lemma~\ref{lemma: entries time is harmonic}, for any $\nu \in \Prob_\mu\ler{Y_Q}$ we have
    \[
        h_\nu := \frac{1}{p_\nu\ler{e}} p_\nu \in \Hcal^{1}(G,\mu).
    \]

    \begin{lemma}\label{lemma: H_mu^1 grow together}
        For every $M > 0$, there exists a finite subset $Q \subset G$ such that for any $h \in \Hcal^{1}(G,\mu)$,
        \[
            \sum_{q \in Q} h\ler{q} > M.
        \]
    \end{lemma}

    \begin{proof}
        First, note the following pointwise bound on every $h \in \Hcal^{1}(G,\mu)$. Since $\mu$ defines a generating random walk, for every $g \in G$ there exists some $k = k\ler{g} \in \NN$ such that $\mu^{*k\ler{g}}\ler{g} > 0$. By $\mu$-harmonicity,
        \[
            1 = h\ler{e} = \sum_{t \in G} \mu^{*k\ler{g}}\ler{t} h\ler{t}
            \geq
            \mu^{*k\ler{g}}\ler{g} h\ler{g},
        \]
        and hence for any $h \in \Hcal^{1}(G,\mu)$ and $g \in G$,
        \begin{equation}\label{eq: bound on H_mu^1}
            h\ler{g} \leq \frac{1}{\mu^{*k\ler{g}}\ler{g}}.
        \end{equation}

        Suppose towards a contradiction that there exists some $M > 0$ such that for every finite subset $Q \subset G$, there exists some $h \in \Hcal^{1}(G,\mu)$ satisfying
        \[
            \sum_{q \in Q} h\ler{q} \leq M.
        \]
        Let $Q_n \nearrow G$ be an exhaustion of $G$ by finite sets, and choose $h_n \in \Hcal^{1}(G,\mu)$ such that
        \[
            \sum_{q \in Q_n} h_n\ler{q} \leq M.
        \]
        By \eqref{eq: bound on H_mu^1}, we may pass to a subsequence $h_{n^\prime}$ that converges pointwise to some function $h \colon G \to [0,\infty)$. Note that for every fixed $m$,
        \[
            \sum_{q \in Q_m} h\ler{q} = \lim_{n^\prime}
            \sum_{q \in Q_m} h_{n^\prime}\ler{q}
            \leq \liminf_{n^\prime}
            \sum_{q \in Q_{n^\prime}} h_{n^\prime}\ler{q} \leq M.
        \]
        Taking $m \to \infty$, we obtain that $h$ is a non-negative function in $\ell^1\ler{G}$. Define $\Tilde{h}\ler{g} := \sum_{t \in G} \mu\ler{t} h\ler{gt}$. Fatou's lemma yields $h \geq \Tilde{h}$ pointwise, while Tonelli's theorem implies that $\| h \|_1 = \| \Tilde{h} \|_1$. Therefore, $h = \Tilde{h}$.

        Thus, $h$ is a non-negative $\mu$-harmonic function in $\ell^1\ler{G}$. Since $G$ is infinite and $\mu$ is generating, $h$ must vanish identically, $h \equiv 0$, which contradicts $h\ler{e} = 1$.
    \end{proof}

    \begin{lemma}\label{lemma: stationary marker lemma}
        Given two finite subsets $I, J \subset G$ and $\epsi > 0$, there exists a finite, containing the unit element, symmetric subset $Q \subset G$ such that $I \subset Q$ and for any $\mu$-stationary measure $\nu \in \Prob_\mu\ler{Y_Q}$,
        \[
            \nu\lers{s \in Y_Q \mid s^{-1}\ler{1} \cap J \ne \emptyset} < \epsi.
        \]
    \end{lemma}

    \begin{proof}
        For each $j \in J$, let $k\ler{j} \in \NN$ be as in \eqref{eq: bound on H_mu^1}. Set
        \[
            B_J := \sum_{j \in J} \frac{1}{\mu^{*k\ler{j}}\ler{j}}
        \]
        and $M := B_J / \epsi$. By Lemma~\ref{lemma: H_mu^1 grow together}, we can choose a finite, containing the unit element, symmetric subset $Q \subset G$ large enough so that $I \subset Q$ and for every $h \in \Hcal^{1}(G,\mu)$,
        \[
            \sum_{q \in Q} h\ler{q} > M.
        \]

        Let $\nu \in \Prob_\mu\ler{Y_Q}$ be a $\mu$-stationary measure, and consider $p_\nu$ and $h_\nu$ as defined above:
        \[
            \begin{split}
                p_\nu\ler{g} &= \nu\lers{s \in Y_Q \mid s\ler{g} = 1}, \\
                h_\nu\ler{g} &= \frac{1}{p_\nu\ler{e}} p_\nu\ler{g} \in \Hcal^{1}(G,\mu).
            \end{split}
        \]
        Note that for every distinct pair $q_1, q_2 \in Q$, the events $\lers{s\ler{q_1} = 1}$ and $\lers{s\ler{q_2} = 1}$ are disjoint, since by symmetry $e \in Qq_1 \cap Qq_2$ and $s^{-1}\ler{1}$ is $Q$-spaced. Therefore, we have
        \[
            1 \geq \nu\lers{s \in Y_Q \mid s^{-1}\ler{1} \cap Q \ne \emptyset}
            = \sum_{q \in Q} p_\nu\ler{q}
            = p_\nu\ler{e} \sum_{q \in Q} h_\nu\ler{q}
            > p_\nu\ler{e} M.
        \]
        Hence,
        \[
            p_\nu\ler{e} < \frac{1}{M} = \frac{\epsi}{B_J}.
        \]
        Finally, combining this inequality with \eqref{eq: bound on H_mu^1} yields the desired bound:
        \[
        \begin{split}
            \nu\lers{s \in Y_Q \mid s^{-1}\ler{1} \cap J \ne \emptyset}
            &\leq \sum_{j \in J} p_\nu\ler{j}
            \\ &= p_\nu\ler{e} \sum_{j \in J} h_\nu\ler{j}
            \\ &\leq p_\nu\ler{e} B_J
            < \epsi.
        \end{split}
        \]
    \end{proof}

    \begin{corollary}
        Theorem~\ref{thm: any RW has universal minimal model} holds.
    \end{corollary}

    \begin{proof}
        Let $K_n \nearrow G$ be an exhaustion of $G$ by finite sets. Using Lemma~\ref{lemma: stationary marker lemma}, we inductively construct three sequences of finite, containing the unit element, symmetric subsets $D_n, E_n, F_n$ exhausting $G$ as described in Subsection~\ref{subsection: the space X_alpha}, with one additional requirement: at step $n$, we apply Lemma~\ref{lemma: stationary marker lemma} with $I = E_n^5$, $J = E_n^4 K_n$, and $\epsi = \frac{1}{2^n}$ to obtain $F_n := Q$ such that for any $\mu$-stationary measure $\nu \in \Prob_\mu\ler{Y_{F_n}}$,
        \begin{equation}\label{eq: key demand stationary}
            \nu\lers{s \in Y_{F_n} \mid s^{-1}\ler{1} \cap E_n^4 K_n \ne \emptyset} < \frac{1}{2^n}.
        \end{equation}

        Now let $\ler{\Omega,\nu}$ be any essentially free $\mu$-stationary Borel $G$-space, and let $U_n \colon \Omega \to Y_{F_n}$ be a weak model map. Since $U_n$ is $G$-equivariant, the push-forward measure ${U_n}_*\nu \in \Prob\ler{Y_{F_n}}$ is also $\mu$-stationary. For any $g \in G$, we eventually have $g \in K_n$, and hence by \eqref{eq: key demand stationary},
        \[
            \sum_{n \in T} {U_n}_*\nu\lers{s \in Y_{F_n} \mid s^{-1}\ler{1} \cap E_n^4 g \ne \emptyset} < \infty.
        \]

        Thus, by Lemma~\ref{lemma: key lemma}, for any sequence $\alpha \in \lers{0,1}^\NN$ with infinitely many $0$'s and $1$'s, the subshift $X_\bullet$ associated with $\bullet = \ler{\textbf{D},\textbf{E},\textbf{F},\alpha}$ is a minimal model for every essentially free $\mu$-stationary Borel $G$-space $\ler{\Omega,\nu}$.
    \end{proof}
    \printbibliography

\end{document}